\documentclass[11pt]{article}
\usepackage[english]{babel}

\usepackage{amsmath}
\usepackage{amssymb}

\newcommand{\R}{   {\ifmmode{{\mathbb R}}\else{$\mathbb R$}\fi}}
\newcommand{\N}{   {\ifmmode{{\mathbb N}}\else{$\mathbb N$}\fi}}

\newcommand{\Lra}{\Longrightarrow}
\newcommand{\lmt}{\longmapsto}

\newcommand{\D}{\mathcal{D}}

\newcommand{\demo}{\noindent\textit{Proof. }}

\newcommand{\sse}{\Longleftrightarrow}

\newcommand{\Li}{\mathcal{L}}
\newcommand{\Lo}{\mathcal{L}_0}
\newcommand{\g}{\mathfrak{g}}
\newcommand{\h}{\mathfrak{h}}
\newcommand{\cc}{\mathcal{C}}
\newcommand{\cci}{\mathcal{C}^{\infty}}
\newcommand{\X}{\mathcal{X}}

\newcommand{\dt}{\frac{d}{dt}}
\newcommand{\dto}{\frac{d}{dt}_{|t=0}}
\newcommand{\ddx}[1]{\frac{\partial}{\partial #1}}
\newcommand{\n}{\mathcal{N}}
\newcommand{\ad}{\mbox{ad}}
\newcommand{\rank}{\mbox{rank}\ }

\begin{document}

\title{Equivalence of Control Systems with Linear Systems on
  Lie Groups and Homogeneous Spaces}
\author{Philippe JOUAN\footnote{LMRS, CNRS UMR 6085, Universit\'e
    de Rouen, avenue de l'universit\'e BP 12, 76801
    Saint-Etienne-du-Rouvray France. E-mail: Philippe.Jouan@univ-rouen.fr}}

\date{\today}

\maketitle

\begin{abstract}
The aim of this paper is to prove that a control affine system on a manifold is equivalent by diffeomorphism
to a linear system on a Lie group or a homogeneous space if and only the vector
fields of the system are complete and generate a finite dimensional
Lie algebra.

A vector field on a connected Lie group is linear if its flow is a one parameter
group of automorphisms. An affine vector field is obtained by adding a
left invariant one. Its projection on a homogeneous space, whenever it exists, is still called affine.

Affine vector fields on homogeneous spaces can be characterized by their Lie brackets with
the projections of right invariant vector fields.

A linear system on a homogeneous space is a system whose drift part is
affine and whose controlled part is invariant.

The main result is based on a general theorem on finite dimensional algebras generated by complete vector fields, closely related to a theorem of Palais, and which have its own interest. The present proof makes use of geometric control theory arguments.

\vskip 0.2cm


\vskip 0.2cm


\end{abstract}

\section{Introduction.}

The aim of this paper is to characterize the class of control systems
\begin{equation}\label{system}
\dot{p}=f(p)+\sum_{j=1}^m u_j g_j(p)
\end{equation}
which are diffeomorphic in the large to a linear system on a Lie group or a homogeneous space. They turn out to be the systems whose vector fields are complete and generate a finite dimensional Lie algebra.

We say that a vector field on a connected Lie group is {\it linear} if its flow is a one parameter group of automorphisms. Linear vector fields on Lie groups are nothing else than the so-called \textit{infinitesimal automorphisms} in the Lie group literature (see for instance \cite{Bourbaki2}). They were first consider in a control theory context by Markus, on matrix Lie groups (see \cite{Markus81}), and then in the general case by Ayala and Tirao (see \cite{AT99}). They are the natural extension to Lie groups of the linear fields on vector spaces and, for this reason, are still called linear.

There is another way, of course equivalent, to define linear vector fields on a connected Lie group $G$. Let us denote by $\g$ its Lie algebra, that is the set of right invariant vector fields. Then a vector field $\X$ is linear if and only if
\begin{equation}
\forall Y\in\g \qquad [\X,Y]\in \g \label{norm}
\end{equation}
and moreover satisfies $\X(e)=0$, where $e$ stands for the identity of $G$. In case where $\X$ satisfies Condition (\ref{norm}) but not $\X(e)=0$ it will be said affine (it is in that case equal to the sum of a linear vector field and a left invariant one). The interest of this second definition is twofold. On the one hand it does not require the knowledge of the flow of $\X$. On the other one it can be extended to homogeneous spaces.

An affine vector field on a homogeneous space $G/H$ is the projection, if it exists, of an affine vector field on $G$. Theorem \ref{homogene} provides a characterization of affine vector fields on homogeneous spaces in terms of Lie brackets with the invariant ones, similar to (\ref{norm}).

A system defined on a homogeneous space,
\begin{equation}
\dot{x}=F(x)+\sum_{j=1}^m u_j Y_j(x), 
\end{equation}
is called linear if the field $F$ is affine and the $Y_j$'s invariant. Linear systems on Lie groups and invariant systems appear as particular cases of this general setting.

The extension of affine vector fields and linear systems to homogeneous spaces is motivated by the fact that the class of systems diffeomorphic to a linear one on a Lie group is rather restrictive, while the class of systems diffeomorphic to linear systems on homogeneous spaces is much wider. More accurately we have (Theorem \ref{maintheorem}, Section \ref{main}):

{\it We assume the family
$\{f,g_1,\dots,g_m\}$ to be transitive.
Then System (\ref{system}) is diffeomorphic to a linear system on a Lie group
or a homogeneous space if and only if the vector fields $f,g_1,\dots,g_m$ are complete and generate a 
finite dimensional Lie algebra.}

In this statement, the transitivity assumption means that there is only one orbit, equal to the state space, under the action of the system. This is not really a limitation since by the Orbit Theorem (see Section \ref{defandnot}) we known that we can always consider the restriction of the system to the orbit through a given point.

The proof makes appeal to Theorem \ref{tecnic} (Section \ref{finite}):
{\it
Let $\Gamma$ be a transitive family of vector fields on a connected  manifold $M$. If all the vector fields belonging to $\Gamma$ are complete,
and if $\Gamma$ generates a finite dimensional Lie algebra $\Li(\Gamma)$,
then $M$ is diffeomorphic to a homogeneous space $G/H$,
where $G$ is a simply connected Lie group whose Lie algebra is isomorphic
to $\Li(\Gamma)$.
By this diffeomorphism $\Li(\Gamma)$ is related to the Lie algebra of invariant
vector fields on $G/H$, and all the vector fields belonging to $\Li(\Gamma)$ are complete.}

Theorem \ref{tecnic} is very closely related to a theorem of Palais (\cite{Palais57}, Theorem III, page 95), but the new proof given here uses control theory ideas like transitivity, rank condition, normal accessibility, Sussmann's Orbit Theorem. An important consequence of Theorem \ref{tecnic} is that the families of vector fields under consideration are Lie determined. 

\vskip 0.2cm

Thus linear systems on Lie groups and homogeneous spaces appear as models for a wide class of systems, and will certainly play an important role in studying topics such as controllability or others. Some results about controllability of linear systems are yet known, see \cite{AT99}, \cite{CM05}, \cite{AS00}.

\vskip 0.2cm

The paper is organized as follows.

In Section \ref{defandnot} the definitions and facts from control theory used in the sequel are recalled.

In Section \ref{Lin} the various definitions of linear and affine vector fields on Lie groups are stated, and their equivalence proved (see Theorem \ref{T1}), as well as their properties, in particular their completeness which is one of the main ingredient used herein. Some of the proofs can be found in the literature (for instance \cite{Bourbaki2}, \cite{AT99}, \cite{CM05}) and are in that case quoted, but to the author's knowledge, Theorem \ref{T1} is nowhere completely proved.

In Section \ref{esp.hom} affine vector fields on homogeneous spaces are introduced, and characterized by their Lie brackets with the projections of right invariant vector fields (Theorem \ref{homogene}).

Section \ref{finite} is devoted to the uppermentioned theorem \ref{tecnic}, and to a corollary (see Corollary \ref{Palais}) where the transitivity assumption is relaxed .

The main result, Theorem \ref{maintheorem}, is stated and proved in Section \ref{main}, which is ended by Corollary \ref{mainbis} where systems diffeomorphic to a linear system on a simply connected Lie group are characterized.

Section \ref{Exemples} begins by some examples of linear vector fields. Then we consider a well known system and show that it is equivalent to a linear one on a homogeneous space of the group Heisenberg. The equivalence is computed.

\vskip 0.2cm

Throughout the paper (in fact from Section \ref{finite}) the vector fields under consideration are only assumed to be $\cc^k$, $k\geq 1$, because the class of differentiability does not matter. The important properties are the completeness of the vector
fields and the fact that they generate a finite dimensional Lie algebra.

\vskip 0.2cm

On the other hand the simply connected spaces will be assumed to be connected. However this definition of simply connectedness, that can be found for instance in \cite{Hochschild65}, is not universal and will be recalled in some statements.


\section{Preliminaries}\label{defandnot}
In this section some standart definitions and facts from control theory are reviewed.

Let $\Gamma=\{g_i;\ i\in I\}$ be a family of $\cc^k$ vector fields on a connected $\cc^{k+1}$ manifold $M$, with $k\geq 1$ (the $g_i$'s are not required to be complete).

Let us denote by $(\gamma^i_t)$ the flow of $g_i$. The \textit{orbit} of $\Gamma$ through a point $p\in M$ is the set of point $q$ for which there exist vector fields $g_{i_1},\dots,g_{i_r}\in\Gamma$, and real numbers $t_1,\dots,t_r$ such that
$$
\gamma^{i_r}_{t_r}\circ\dots\circ\gamma^{i_1}_{t_1}(p)
$$
is defined and equal to $q$. Let us recall the Sussmann's Theorem:

\vskip 0.3cm

\noindent{\bf Orbit Theorem}.
The orbit of $\Gamma$ through each point $p$ of $M$ is a connected submanifold of $M$.

\vskip 0.3cm

In the original proof of Sussmann the vector fields are assumed to be $\cci$ (see \cite{Sussmann73}), but a proof for $\cc^k$ vector fields, with $k\geq 1$, can be found in \cite{Jurdjevic97}.

\vskip 0.3cm

This family of vector fields is said to be \textit{transitive} if the orbit through each point $p$ of $M$ is equal to $M$, that is if $M$ is the only orbit of $\Gamma$. 

\vskip 0.3cm

Let $V^k(M)$ stand for the space of $\cc^k$ vector fields on $M$. It is not a Lie algebra whenever $k<+\infty$.
But it may happen that all the Lie brackets of elements of $\Gamma$ of all finite lengths exist and are also $\cc^k$. In that case the subspace of $V^k(M)$ spanned by these Lie brackets is a Lie algebra and we will say that \textit{the family $\Gamma$ generates a Lie algebra}. This last will be denoted by $\Li(\Gamma)$.

\vskip 0.3cm

Let us assume that $\Gamma$ generates a Lie algebra, and let us consider the \textit{rank} of $\Gamma$ at each point $p\in M$, that is the dimension of the subspace of $T_pM$, the tangent space to $M$ at $p$, spanned by the vectors $\gamma(p),\ \ \gamma\in \Li(\Gamma)$. The so-called \textit{rank condition} asserts that the family $\Gamma$ is transitive as soon as its rank is maximum, hence equal to $\dim M$, at each point.

\vskip 0.3cm

To finish let us recall the definition of Lie-determined systems (see for instance \cite{Jurdjevic97}): the family $\Gamma$ is said to be \textit{Lie-determined} if at each point $p\in M$, the rank of $\Gamma$ at $p$ is equal to the dimension of the orbit of $\Gamma$ through $p$.


\section{Linear and affine vector fields on Lie groups.}\label{Lin}

Let $G$ be a {\bf connected} Lie group, and $\g$ its Lie algebra, that is the set of
right invariant vector fields. Let us denote by
$V^{\omega}(G)$ the set of analytic vector fields on $G$. The
normalizor of $\g$ in $V^{\omega}(G)$ is by definition
$$
\n=\mbox{norm}_{V^{\omega}(G)}\g=\{F\in V^{\omega}(G)/\ \forall Y\in \g \qquad
[F,Y]\in\g\}
$$

\newtheorem{D1}{Definition}
\begin{D1}
A vector field $F$ on $G$ is said to be {\bf affine} if it belongs to $\n$.

Such a vector field $F$ is said to be linear if it moreover verifies $F(e)=0$, where $e$ stands for the identity of $G$.
\end{D1}

In other words the restriction of $ad(F)$ to $\g$, also denoted by $ad(F)$, is a derivation of
$\g$.
From the Jacobi identity, it is clear that $\n$ is a Lie subalgebra
of $V^{\omega}(G)$, and that the mapping $F\lmt ad(F)$ is a Lie
algebra morphism from $\n$ into $\D(\g)$, the set of derivations of
$\g$.

\newtheorem{P1}{Proposition}
\begin{P1}\label{noyau}
(See \cite{AT99}) The kernel of the mapping $F\lmt ad(F)$ is the set of left invariant
vector fields. An affine vector field $F$ can be uniquely decomposed
into a sum
$$
F=\X+Z
$$
where $\X$ is linear and $Z$ left invariant.
\end{P1}

This proposition is no longer true whenever the group $G$ is not connected.

\demo

Let $Z$ be an affine vector field whose bracket with any element of
$\g$ vanishes. Its flow $z_t$ commutes with the one of any $Y\in
\g$. This writes
$$\forall x\in G \qquad z_t(\exp\ (sY)x)=\exp\ (sY)z_t(x) $$
for all $t,s\in \R$ for which this makes sense. Fix $x\in G$. There
exist $Y_1\dots Y_k\in \g$ such that
$$x=\exp\ (Y_1)\dots \exp\ (Y_k),$$
thanks to the connectedness of $G$. Therefore
$$z_t(x)=z_t(\exp\ (Y_1)\dots \exp\ (Y_k))=\exp\ (Y_1)\cdots \exp\ (Y_k)z_t(e)=xz_t(e)$$
for $t$ sufficiently small, and
$$
\begin{array}{ll}
Z_x & =\dto z_t(x)\\
           & =\dto xz_t(e)\\
           & =T_eL_x.Z_e
\end{array}
$$
where $T_eL_x$ stands for the tangent mapping at the identity $e$ of
the left translation $L_x$. This proves that $Z$ is left invariant and, the converse being obvious, the first part of the proposition.

 For the second one let $Z$ be the left invariant vector field defined by $Z_e=F(e)$. Then the vector field $\X=F-Z$ is clearly linear.

\hfill $\Box$

\newtheorem{T1}{Theorem}
\begin{T1}\label{T1}
Let $\X$ be a vector field on a connected Lie group $G$. The following
conditions are equivalent:
\begin{enumerate}
\item
$\X$ is linear;
\item
the flow of $\X$ is a one parameter group of automorphisms of $G$;
\item
$\X$ verifies
\begin{equation}
 \forall x,x'\in G \qquad \X_{xx'}=TL_x.\X_{x'}+TR_{x'}.\X_x \label{Bourbak}
\end{equation}
\end{enumerate}
The second item implies that a linear vector field on a connected Lie group is complete.
\end{T1}

\demo
\begin{enumerate}
 \item[1$\ \Rightarrow$ 3]
The vector field $\X$ is linear, hence for all $Y\in \g$,
the Lie bracket $[\X,Y]$ is right invariant. Therefore we have
$$
\forall x\in G \quad
[\X,Y]=R_{x*}[\X,Y]=[R_{x*}\X,R_{x*}Y]=[R_{x*}\X,Y].$$

This proves firstly that $[R_{x*}\X,Y]$ is right invariant for all $Y\in\g$, hence that the vector field $R_{x*}\X$ is affine, and secondly that $\ad(R_{x*}\X)=\ad(\X)$ hence following Proposition \ref{noyau} that
\begin{equation}
 \label{RG}
R_{x*}\X=\X+Z
\end{equation}
where $Z$ is left invariant. This last is characterized by
$$Z_e=(R_{x*}\X)_e=TR_x.\X_{x^{-1}}$$
hence for all $x'\in G$
$$
Z_{x'}=TL_{x'}\ TR_x.\X_{x^{-1}}=TR_x\ TL_{x'}.\X_{x^{-1}}.
$$
Considering  $(R_{x*}\X)_{x'}=TR_x.\X_{x'x^{-1}}$, Equality \ref{RG} evaluated at the point $x'$ becomes
$$
TR_x.\X_{x'x^{-1}}=\X_{x'}+TR_x\ TL_{x'}.\X_{x^{-1}}.
$$
It remains to apply $TR_{x^{-1}}$ to obtain
$$
\X_{x'x^{-1}}=TR_{x^{-1}}.\X_{x'}+TL_{x'}.\X_{x^{-1}}
$$
and to replace $x^{-1}$ by $x$ to obtain Equality \ref{Bourbak}.

\item[3$\ \Rightarrow$ 2] (see also \cite{Bourbaki2})
Let us denote by $\varphi_t$ the flow of $\X$, defined on a domain of $\R\times G$. The curve
$$t\lmt \varphi_t(x)\varphi_t(x')$$
is defined on an open interval containing $0$, and takes the value $xx'$ at $t=0$. Moreover 
$$
\begin{array}{ll}
\dt\varphi_t(x)\varphi_t(x') & = TL_{\varphi_t(x)}.\X_{\varphi_t(x')}+TR_{\varphi_t(x')}.\X_{\varphi_t(x)}\\
                             & = \X_{\varphi_t(x)\varphi_t(x')}
\end{array}
$$
This proves that 
$$\varphi_t(x)\varphi_t(x')=\varphi_t(xx')$$
as soon as the left-hand side exists. It remains to show that $\X$ is complete. To begin with, notice that Equality \ref{Bourbak}, evaluated at $x=x'=0$, implies $\X_e=0$.

Let $x\in G$ and $t\in\R$. Since $\X_e$ vanishes, $\varphi_t$ is defined on an open neighborhood $V_t$ of $e$. The group $G$ being connected, it is generated by this neighborhood. Therefore there exist
$x_1,\dots ,x_n\in V_t$ such that $x=x_1\dots x_n$, and 
$$
\varphi_t(x)=\varphi_t(x_1\dots x_n)=\varphi_t(x_1)\dots\varphi_t(x_n)
$$
is well defined. This proves that $\X$ is complete and that $\varphi_t$ is an automorphism of $G$ for all $t\in\R$.
\item[2$\ \Rightarrow$ 1]
Let $(\varphi_t,\ t\in \R)$ be a one parameter group of automorphisms of $G$, and $\X$ its infinitesimal generator.
For all right invariant vector fields $Y$, we have
\begin{equation}
[\X,Y]_e=\dto T_{\varphi_t(e)}\varphi_{-t}.Y_{\varphi_t(e)}=\dto T_e\varphi_{-t}.Y_e \label{tzero}
\end{equation}
since $\varphi_t(e)=e$ for all $t\in \R$. Considering $\varphi_{-t}\circ R_{\varphi_t(x)}=R_x\circ\varphi_{-t}$ we have at any point $x$:
$$
\begin{array}{ll}
[\X,Y]_x & =\dto T_{\varphi_t(x)}\varphi_{-t}.Y_{\varphi_t(x)}\\
           & =\dto T_{\varphi_t(x)}\varphi_{-t}\ T_eR_{\varphi_t(x)}.Y_e\\
           & =\dto T_eR_x\ T_e\varphi_{-t}.Y_e\\
           & =T_eR_x.[\X,Y]_e,
\end{array}
$$
The vector field $\X$ is therefore affine, and consequently linear since $\varphi_t(e)=e$ for all $t\in\R$.

\end{enumerate}
\hfill $\Box$

\noindent{\bf Notations}. Here and subsequently the flow of a linear
vector field $\X$ will be denoted by $(\varphi_t)_{t\in\R}$.

\vskip 0.1cm

To a given linear vector field $\X$, one can associate the derivation $D$ of $\g$ defined by:
$$
\forall Y\in \g \qquad DY=-[\X,Y],
$$
that is $D=-\ad(\X)$.
The minus sign in this definition comes from the formula
$[Ax,b]=-Ab$ in $\R^n$. It also enables to avoid a minus sign in the equality
$$
\forall Y\in\g \quad \forall t\in \R \qquad \varphi_t(\exp
Y)=\exp(e^{tD}Y),
$$
stated in the forthcoming proposition \ref{CVL 1}.

\vskip 0.2cm

\noindent{\bf Example: the inner derivations}

Let $X\in \g$ be a right invariant vector field. We denote by $\mathcal{I}$ the diffeomorphism of $G$ defined by $x\lmt x^{-1}$.
The vector field $\mathcal{I}_*X$ is left invariant and equal to 
$-X_e$ at $e$. Therefore the vector field
$$
\X=X+\mathcal{I}_*X
$$
is linear. Indeed $\X$ belongs to $\n$, because for all $Y\in \g$, we have $[\X,Y]=[X+\mathcal{I}_*X,Y]=[X,Y]\in\g$, and satisfies moreover $\X(e)=0$.
The derivation associated to $\X$ is inner since it is equal to $-\ad(\X)=-\ad(X)$.

\vskip 0.2cm
This shows also that given an \textit{inner} derivation $D=-\ad(X)$, there always exists a linear vector field on $G$ whose associated derivation is $D$. This is no longer true in the general case, but remains true whenever $G$ is simply connected. This fact is crucial in the sequel.

\newtheorem{simplyconnected}[T1]{Theorem}
\begin{simplyconnected}\label{simplyconnected}
The group $G$ is assumed to be (connected and) simply connected. Let $D$ be a derivation of its Lie algebra $\g$. Then there exists one and only one  linear vector field on $G$ whose associated derivation is $D$.
\end{simplyconnected}

\demo
The proof is essentially contained in \cite{Bourbaki2}, lemme 4, page 250.

\hfill $\Box$

\noindent{\bf Important remark.} Under the assumption that the group $G$ is connected, it was stated in Proposition \ref{noyau} that an affine vector field $F$ can be decomposed into $F=\X+Z$, with $\X$ linear and $Z$ left invariant. This decomposition is natural because $\ad(\X)=\ad(F)$, but it may be useful to decompose $F$ into a linear part and a right invariant one. Since $\mathcal{I}_*Z$ is right invariant and $Z+\mathcal{I}_*Z$ linear, we can write
$$
F=\widetilde{\X}-\mathcal{I}_*Z \quad\mbox{ with }\widetilde{\X}=\X+Z+\mathcal{I}_*Z.
$$
Of course $\ad(\widetilde{\X})=\ad(\X)+\ad(\mathcal{I}_*Z)\neq \ad(\X)=\ad(F)$.\\
This remark will be used in Section \ref{esp.hom}.

\newtheorem{CVL 1}[P1]{Proposition}
\begin{CVL 1}\label{CVL 1}
For all $t\in\R$
$$T_e\varphi_t=e^{tD}$$
and consequently
$$\forall Y\in\g \quad \forall t\in \R \qquad  \varphi_t(\exp Y)=\exp (e^{tD}Y).$$
\end{CVL 1}

\demo
\begin{enumerate}
\item
Let us first prove the equality
$$\dt T_e\varphi_t.Y_e=DT_e\varphi_t.Y_e$$
This equality has been previously stated at $t=0$ (see Equality (\ref{tzero}) in the proof of Theorem \ref{T1}).
In general
$$
\begin{array}{ll}
\dt T\varphi_t.Y_e & =\frac{d}{ds}_{|s=0}T_e\varphi_{t+s}.Y_e\\
                      & =\frac{d}{ds}_{|s=0}T_e\varphi_s.T_e\varphi_t.Y_e\\
                      & =DT_e\varphi_t.Y_e
\end{array}
$$
\item
From the equality proved above, the first formula of the proposition
is immediate. For the second one, remark that $\varphi_t$ is a Lie
group morphism. Therefore
$$
\begin{array}{ll}
\varphi_t(\exp Y) & =\exp (T_e\varphi_t.Y)\\
                  & =\exp (e^{tD}Y)
\end{array}
$$
\end{enumerate}

\hfill $\Box$

To finish this section notice the following proposition:

\newtheorem{Complete}[P1]{Proposition}
\begin{Complete}
An affine vector field on a connected Lie group is complete.
\end{Complete}

\demo
This proposition is a consequence of the forthcoming Theorem
\ref{tecnic}, but it can be proved in an elementary way. Indeed let $F$ be
an affine vector field and $\X+Z$ its decomposition into a linear
vector field $\X$ and a left invariant one $Z$. Let us
denote by $t\lmt e(t)$ the maximal trajectory of $F$ through the identity $e$,
defined on an interval $]a,b[$. One can verify, using the third
characterization of linear vector fields (see Theorem \ref{T1}), that $t\lmt \varphi_t(x)e(t)$
is the trajectory of $F$ through the point $x\in G$, also defined on
$]a,b[$. Let us assume $b<+\infty$ and let us choose $x=e(b/2)$. Then
$$
t\lmt\varphi_{t-\frac{b}{2}}(e(\frac{b}{2}))e(t-\frac{b}{2})
$$
is the trajectory of $F$ through $x=e(b/2)$ at $t=b/2$. Therefore the
trajectory $t\lmt e(t)$ can be extended up to $3b/2$, a contradiction.

\hfill $\Box$


\section{Affine vector fields on homogeneous spaces.}\label{esp.hom}

Let $H$ be a closed subgroup of $G$. The homogeneous space $G/H$ is the manifold of left cosets of $H$, and we denote by $\Pi$ the projection of $G$ onto $G/H$. For any right invariant vector field $Y\in\g$, the projection $\Pi_*Y$ of $Y$ onto $G/H$ is always well defined, and will be refered to as an invariant vector field on $G/H$.
It is well known that the set of such vector fields, $\Pi_*\g=\{\Pi_*Y;\ Y\in \g\}$, is
a Lie algebra and that $\Pi_*$ is a Lie algebra morphism from $\g$ onto
$\Pi_*\g$.

Let $\X$ be a linear vector field on
$G$. We investigate the existence on $G/H$ of a vector field $\Pi$-related to $\X$. Such a vector field exists if and only if
$$
\forall x\in G,\ \forall y\in H,\ \forall t\in\R \qquad \Pi(\varphi_t(xy))=\Pi(\varphi_t(x)).
$$
But $\Pi(\varphi_t(xy))=\varphi_t(x)\varphi_t(y)H$, and the preceding condition is equivalent to
$$
\forall y\in H,\ \forall t\in\R \qquad \varphi_t(y)\in H.
$$
Thus $\X$ is $\Pi$-related to a vector field on $G/H$ if and only if $H$ is invariant under the flow of $\X$, therefore if and only if $\X$ is tangent to $H$.

In the particular case where $H$ is a discrete subgroup of $G$, this amounts to the condition that $\X$ vanishes everywhere on $H$, or that $H$ is included in the set of fixed points of $\X$.

Assume now $H$ to be connected, and denote by $\h$ its Lie algebra. Since the elements of $H$ are products of exponentials the invariance of $H$ under $\X$ writes
$$
\forall Y\in \h,\ \forall t\in\R  \qquad \varphi_t(\exp Y)=\exp(e^{tD}Y)\in H.
$$
This is equivalent to $\forall Y\in \h,\ \forall t\in\R$, $e^{tD}Y\in\h$, and finally to the invariance of $\h$ under $D$.

\newtheorem{Projecthomogene}[P1]{Proposition}
\begin{Projecthomogene}\label{Projecthomogene}
Let $H$ be a closed subgroup of $G$, $G/H$ the homogeneous space of left cosets of $H$, and $\Pi$ the projection of $G$ onto $G/H$.

A linear vector field $\X$ on
$G$ is $\Pi$-related  to a vector field on $G/H$ if and only if $H$ is invariant under $\X$.

If $H$ is discrete, this condition holds if and only if $H$ is included in the set of fixed points of $\X$.

If $H$ is connected it is equivalent to the invariance of its Lie algebra $\h$ under the derivation $D$ associated to $\X$.
\end{Projecthomogene}

Under these conditions, the projection of $\X$ onto $G/H$ will be denoted by $\Pi_*\X$.

Let us now consider an affine vector field $F$ on $G$. It is equal to 
$$
F=\X+Y
$$
where $\X$ is linear and $Y$ right invariant. This decomposition (see ``important remark'' in Section \ref{Lin}) is chosen in order to ensure that the projection $\Pi_*Y$
of $Y$ onto $G/H$ is well defined. Then $F$ is $\Pi$-related to a vector
field on $G/H$ if and only if $\Pi_*\X$ exists. In that case
$\Pi_*F=\Pi_*\X+\Pi_*Y$ will stand for the projection of $F$ onto
$G/H$.

\newtheorem{Affinehomogene}[P1]{Proposition}
\begin{Affinehomogene}\label{Affinehomogene}
Let $H$ be a closed subgroup of $G$, $G/H$ the homogeneous space of left cosets of $H$, and $\Pi$ the projection of $G$ onto $G/H$.

Let $F$ be an affine vector field on $G$ and $F=\X+Y$ its
decomposition into a linear vector field $\X$ and a right invariant one $Y$.

Then $F$ is $\Pi$-related  to a vector field on $G/H$ if and only if
this holds for $\X$, hence if and only if $H$ is invariant under $\X$.
\end{Affinehomogene}

The next task is to define and characterize affine vector fields on connected homogeneous
spaces. In general a homogeneous space is defined as a manifold on which a Lie group acts smoothly and transitively. For our purpose it is more convenient to use the following equivalent definition: a (connected) homogeneous space $M$ is a manifold diffeomorphic to a quotient $G/H$, where $G$ is a (connected) Lie group and $H$ a closed subgoup of $G$.

There are two remarks to make about the choice of the Lie group $G$.
\begin{enumerate}
 \item
We can assume $G$ to be simply connected. If not let $\widetilde{G}$
be the universal covering of $G$, $\sigma$ the projection of
$\widetilde{G}$ onto $G$ and
$\widetilde{H}=\sigma^{-1}(H)$. Then $G/H$ is
diffeomorphic to $\widetilde{G}/\widetilde{H}$.
\item
We can assume that $\dim \Pi_*\g=\dim \g$. If not let $\mathfrak{k}$
be the kernel of $\Pi_*$ and let $K$ be the connected Lie subgroup of $G$
whose Lie algebra is $\mathfrak{k}$. The subgroup $K$ is normal and
included in $H$ because $\mathfrak{k}$ is an ideal of $\g$ included in
the Lie algebra $\h$ of $H$ (it is easy to see that $Y\in
\mathfrak{k}$ if and only if $\forall x\in G$, $\forall t\in \R$,
$x^{-1}\exp(tY)x\in H$). Now $G/K$ is a simply connected (because $K$ is connected) Lie group and
$$
G/H\sim (G/K)/(H/K).
$$
\end{enumerate}
We can therefore restrict ourselves to homogeneous spaces $G/H$ where $G$ is simply connected and $H$ is a closed subgroup of $G$ such that
$$
\dim \Pi_*\g=\dim \g
$$

\newtheorem{VFaffine}[D1]{Definition}
\begin{VFaffine}\label{VFaffine}
Let $G$ be a simply connected Lie group, $H$ a closed subgroup of $G$, such that $\dim \Pi_*\g=\dim \g$, where
$\Pi$ stands for the projection of $G$ onto $G/H$.
A vector field $f$ on $G/H$ is said to be affine if it is $\Pi$-related
to an affine vector field of $G$.
\end{VFaffine}

It is clear that $\Pi_*\g$ is invariant for the Lie bracket
with any affine vector field $\Pi_*F$. Let us now state and prove the
converse statement, and thus characterize the affine vector fields on
homogeneous spaces.

\newtheorem{homogene}[T1]{Theorem}
\begin{homogene}\label{homogene}
Let $G$ be a (connected and) simply connected Lie group, $H$ a closed subgroup of $G$, such that $\dim \Pi_*\g=\dim \g$, where
$\Pi$ stands for the projection of $G$ onto $G/H$.

A vector field $f$ on $G/H$ is affine if and only if
$$
\forall Y\in \g \qquad [f,\Pi_*Y]\in\Pi_*\g,
$$
that is if and only if $\Pi_*\g$ is $\ad_f$-invariant.
\end{homogene}

\demo The necessary part being clear, let us prove the converse and,
for this purpose, let us begin by some preliminary remarks.
\begin{enumerate}
\item
We can assume without loss of generality $f(H)=0$ since we can add to $f$ an invariant vector
field $\Pi_*Y$ that verifies $\Pi_*Y(H)=-f(H)$.
\item
Let us notice that an invariant vector field $Y$ belongs to the Lie
algebra $\h$ of $H$ if and only if $\Pi_*Y(H)=0$. Indeed
$$
\begin{array}{ll}
Y\in\h  &  \sse \forall t\in \R \quad \exp(tY)\in H\\
        &  \sse \forall t\in \R \quad \exp(tY)H=H\\
        &  \sse \Pi_*Y(H)=0.
\end{array}
$$
\end{enumerate}

By assumption $f$ induces a derivation $\widetilde{D}$ on $\Pi_*\g$,
defined by $\widetilde{D}=-\ad_f$. Since $\g$ and $\Pi_*\g$ are isomorphic, we can
define the derivation $D$ on $\g$ by the equality $\Pi_*\circ
D=\widetilde{D}\circ\Pi_*$. Moreover $G$ being simply connected there
exists a (unique) linear vector field $\X$ on $G$ associated to $D$.

Let $Y\in\h$. We have
$$
\Pi_*[Y,\X]=\Pi_*(DY)=\widetilde{D}(\Pi_*Y)=[\Pi_*Y,f].
$$
But $\Pi_*Y(H)=0$ because $Y\in \h$, and by assumption
$f(H)=0$. Therefore $[\Pi_*Y,f](H)=0$, and $[Y,\X]$ belongs to $\h$.

This proves that $\h$ is invariant under $D$. Let $H_0$ denote the
connected component of $e$ in $H$. Following Proposition
\ref{Projecthomogene}, the linear vector field $\X$ is $\Pi$-related
to a vector field $\widetilde{f}$ on $G/H_0$.

If $H_0=H$, then $\widetilde{f}=f$. Indeed $\widetilde{f}$ and $f$
verify
$$
\begin{array}{l}
(i)\ \ \forall g\in \Pi_*\g \qquad [\widetilde{f},g]=[f,g]\\
(ii)\ \widetilde{f}(H)=f(H)=0,
\end{array}
$$
and according to the forthcoming Lemma \ref{transitive} this implies the equality of
$\widetilde{f}$ and $f$.

If $H$ is not connected then $G/H_0$ is a covering space of $G/H$, and
$f$ can be lift to a vector field $f'$ on $G/H_0$. By the previous
method we obtain $\widetilde{f}=f'$. Moreover $f$ is related to $f'$
by the projection of $G/H_0$ onto $G/H$, and therefore to $\X$. More
accurately the equality between $\widetilde{f}$ and $f'$ implies that
$\widetilde{f}$ vanishes on $H/H_0$, because so does $f'$. Therefore the
connected components of $H$ are invariant under $\X$, and according to Proposition
\ref{Projecthomogene}, $\X$ is $\Pi$-related to a vector field on
$G/H$ which is nothing else than $f$.

\hfill $\Box$

\newtheorem{transitive}{Lemma}
\begin{transitive} \label{transitive}
Let $(g_i)_{i\in I}$ be a transitive family of vector fields (see Section \ref{defandnot}) on a
connected manifold $M$, and let $f$ be a vector field on $M$ that
satisfies
$$
\begin{array}{l}
(i)\ \ \forall i\in I \qquad [f,g_i]=0\\
(ii)\ \exists x_0\in M \qquad f(x_0)=0.
\end{array}
$$
Then $f=0$.
\end{transitive}

\demo Let $\varphi^i_t$ denote the flow of $g_i$. Let $x\in M$. By
assumption there exist $i_1,\dots, i_r\in I$ and $t_1,\dots,t_r\in \R$
such that
$$
x=\varphi^{i_r}_{t_r}\circ\dots\circ\varphi^{i_1}_{t_1}(x_0).
$$
The flow $\psi_t$ of $f$ commutes with $\varphi^i_t$, $\forall i\in I$.
Hence for all $t$ sufficiently small
$$
\begin{array}{ll}
\psi_t(x)  &
=\varphi^{i_r}_{t_r}\circ\dots\circ\varphi^{i_1}_{t_1}(\psi_t(x_0))\\
  & 
=\varphi^{i_r}_{t_r}\circ\dots\circ\varphi^{i_1}_{t_1}(x_0)\\
  & =x.
\end{array}
$$
Therefore $f(x)=0$.

\hfill $\Box$

Next Proposition is obvious but useful in the sequel.

\newtheorem{Comp}[P1]{Proposition}
\begin{Comp}\label{Comp}
 An affine vector field on a homogeneous space is complete.
\end{Comp}

We can now state the definition of general linear systems. They are the systems
\begin{equation}
 \dot{x}=F(x)+\sum_{j=1}^m u_j Y_j(x)  \label{linear system}
\end{equation}
on homogeneous spaces $G/H$, where the field $F$ is affine and the $Y_j$'s invariant. Linear systems on Lie groups, obtained when the subgroup $H$ is normal, and invariant systems, obtained when the vector field $F$ is invariant, are two particular cases of this general setting.


\section{Finite dimensional algebras of vector fields}\label{finite}

Let $\Gamma=\{g_i;\ i\in I\}$ be a family of vector fields on a connected manifold $M$. All the vector fields $g_i$ belonging to $\Gamma$ are  assumed to be $\cc^k$, for a common $k\geq 1$, and $M$ is therefore at least $\cc^{k+1}$.

Recall that the family $\Gamma$ is said to generate a Lie algebra if all the Lie brackets of elements of $\Gamma$ of all finite lengths exist and are also $\cc^k$, and that we define in that case the Lie algebra $\Li(\Gamma)$ as the subspace of $V^k(M)$ spanned by these Lie brackets.

\newtheorem{tecnic}[T1]{Theorem}
\begin{tecnic} \label{tecnic}
Let $\Gamma$ be a family $\cc^k$ vector fields on a connected  manifold $M$. If
\begin{enumerate}
\item[(i)] all the vector fields belonging to $\Gamma$ are complete,
\item[(ii)] $\Gamma$ generates a finite dimensional Lie algebra $\Li(\Gamma)$,
\item[(iii)] the family $\Gamma$ is transitive,
\end{enumerate}
then $M$ is $\cc^{k+1}$ diffeomorphic to a homogeneous space $G/H$,
where $G$ is a (connected and) simply connected Lie Group whose Lie algebra is isomorphic
to $\Li(\Gamma)$, and $H$ is a closed subgroup of $G$.

By this diffeomorphism $\Gamma$ is related to a set of invariant vector
fields, and $\Li(\Gamma)$ to the Lie algebra of invariant
vector fields on $G/H$.

Moreover all the vector fields belonging to $\Li(\Gamma)$ are complete.
\end{tecnic}

\demo

\begin{enumerate}
\item
Let $G$ be a (connected and) simply connected Lie group whose Lie
algebra $\g$ is isomorphic to $\Li(\Gamma)$, $\Li$ in short, and let us denote by
$$
L:\quad \Li\lmt \g
$$
this Lie algebra isomorphism.

\item
In the product $M\times G$ consider the distribution spanned by the
family of vector fields $\{(g, L(g));\ \ g\in
\Li\}$. This distribution is involutive, and its rank is constant,
equal to $\dim(G)=\dim(\Li)$. Henceforth it is completely integrable (the proof of the Frobenius theorem for $\cc^k$ vector fields on a $\cc^{k+1}$ manifold with $k\geq 1$ can be found for instance in \cite{Malliavin72}).

Let us fix an arbitrary point $p_0$ in $M$, and let $S$ be the leaf of the foliation through the point $(p_0,e)$
(where $e$ stands for the identity element of $G$). We denote by $\Pi_1$
(resp. $\Pi_2$) the projection of $S$ onto $M$ (resp. onto $G$). We are going to prove that $\Pi_2$ is a diffeomorphism.

\item 
\textit{Notations}. For $g_i\in \Gamma=\{g_i;\ i\in I\}$ we denote by $Y_i=L(g_i)$ the corresponding vector field on $G$. The flows of $g_i$ and $Y_i$ are respectively denoted by 
$$
(t,p)\lmt \gamma^i_t(p) \ \mbox{ and }\ 
(t,x)\lmt \exp(t Y_i)x.
$$

\item First of all let us show that $\Pi_2$ is onto. Since $L$ is a Lie algebra isomorphism and $\Gamma$ generates $\Li$, the Lie algebra $\g$ of $G$ is generated by the family $L(\Gamma)=\{Y_i;\ i\in I\}$. This family is therefore transitive on $G$ which is connected. Let $x\in G$. There exists an integer $r$, indices $i_1,\dots,i_r\in I$ and real numbers $t_1,\dots,t_r$ such that
$$
x=\exp(t_rY_{i_r})\dots\exp(t_1Y_{i_1}).
$$
Thanks to the completeness assumption of the elements of $\Gamma$ the point
$$
p=\gamma_{t_r}^{i_r}\circ\dots\circ\gamma_{t_1}^{i_1}(p_0)
$$
is well defined. Moreover the vector fields $(g_i,Y_i)$ are tangent to $S$ and the point $(p,x)$ belongs to $S$. Its projection onto $G$ is $x$, and this proves the surjectivity of $\Pi_2$.

\item
Let us now prove that $\Pi_2$ is a covering map. Since the family $L(\Gamma)=\{Y_i;\ i\in I\}$ generates the Lie algebra $\g$ of $G$, the identity $e$ is normally accessible from $e$ (for the family $\{\pm Y_i;\ i\in I\}$, see \cite{Jurdjevic97} for the notion of normal accessibility, in particular  Corollary 1 page 154). Let $n$ be the dimension of $G$. We can find indices $i_1,\dots,i_n\in I$ and real numbers $t_1,\dots,t_n>0$ such that the mapping
$$
(s_1,\dots,s_n)\lmt \exp((t_n+s_n)Y_{i_n})\dots\exp((t_1+s_1)Y_{i_1})
$$
is a local diffeomorphism at $(0,\dots,0)\in \R^n$.

Let $x_0=\exp(-t_1Y_{i_1})\dots\exp(-t_nY_{i_n})$. Then the mapping $\Psi$ defined by
$$
(s_1,\dots,s_n)\lmt \exp((t_n+s_n)Y_{i_n})\dots\exp((t_1+s_1)Y_{i_1})x_0
$$
is also a local diffeomorphism, and satisfies $\Psi(0,\dots,0)=e$.

We can choose a neighbourhood of $0$ in $\R^n$ sent diffeomorphically by $\Psi$ onto an open and connected neighbourhood $V$ of $e$ in $G$.
Let us denote by $\tilde{\Psi}$ the similar mapping from $\R^n\times M$ into $M$, that is
$$
(s_1,\dots,s_n)\lmt \gamma^{i_n}_{t_n+s_n}\circ\dots\circ\gamma^{i_1}_{t_1+s_1}\circ\gamma^{i_1}_{-t_1}\circ\dots\circ\gamma^{i_n}_{-t_n}
$$

Let $x$ be a given point in $G$, and for every $p\in M$ such that $(p,x)\in S$ let $\sigma_p$ be the mapping
$$
\begin{array}{cccc}
\sigma_p: &  Vx & \lmt & S\\
          &  y  & \lmt & (\tilde{\Psi}(\tau)(p),y)
\end{array}
$$
where $\tau=\Psi^{-1}(yx^{-1})$. The neighbourhood $Vx$ of $x$ is evenly covered by $\{\sigma_p(Vx);\ p\in M \mbox{ and } (p,x)\in S\}$. Indeed it is clear that
$\Pi_2\circ\sigma_p$ is the identity of $Vx$ and that the sets $\sigma_p(Vx)$ cover $\Pi_2^{-1}(Vx)$. Let us show that they are mutually disjoint. If not we can find $\tau$, $\tau'$, $p$, $p'$, such that $\sigma_p(\Psi(\tau).x)=\sigma_{p'}(\Psi(\tau').x)$. But
$$
\begin{array}{ll}
\sigma_p(\Psi(\tau).x)=\sigma_{p'}(\Psi(\tau').x) & \Longleftrightarrow (\tilde{\Psi}(\tau)(p),\Psi(\tau).x)\\
                    & \qquad \quad\qquad\qquad=(\tilde{\Psi}(\tau')(p'),\Psi(\tau').x)\\
                    & \Longrightarrow \Psi(\tau).x=\Psi(\tau').x\\
                    & \Longrightarrow \tau=\tau'\\
                    & \Longrightarrow \tilde{\Psi}(\tau)(p)=\tilde{\Psi}(\tau)(p')\\
                    & \Longrightarrow p=p'.

\end{array}
$$
This proves that $\Pi_2$ is a covering of $G$ by $S$ (notice that $S$ is connected and locally connected). But $G$ is simply connected and $\Pi_2$ is therefore a diffeomorphism. In particular, given a point $x\in G$, there is one and only one point $p\in M$ for which $(p,x)\in S$.

\item

The next task is to prove that the left translations of $G$ induce a group action on $M$. Let $y\in G$ and $(p,x)\in S$. There exists an unique point $q\in M$ such that $(q,yx)$ belongs to $S$. Let us show that $q$ depends only on $p$ and $y$ but not on a particular choice of $x$. Thanks again to the transitivity of the family $L(\Gamma)$, there exists an integer $r$, indices $i_1,\dots,i_r\in I$ and real numbers $t_1,\dots,t_r$ such that
$$
y=\exp(t_rY_r)\dots\exp(t_1Y_1).
$$
Then the equality $q=\gamma_{t_r}^{i_r}\circ\dots\circ\gamma_{t_1}^{i_1}(p)$ holds.

Indeed $(\gamma_{t_r}^{i_r}\circ\dots\circ\gamma_{t_1}^{i_1}(p),\exp(t_rY_r)\dots\exp(t_1Y_1)x)$ belongs to $S$, and $\gamma_{t_r}^{i_r}\circ\dots\circ\gamma_{t_1}^{i_1}(p)$ is therefore the only point $q$ such that \\ $(q,\exp(t_rY_r)\dots\exp(t_1Y_1)x)=(q,yx)$ belongs to $S$.

We will denote by $\rho_y$ the diffeomorphism $\gamma_{t_r}^{i_r}\circ\dots\circ\gamma_{t_1}^{i_1}$ of $M$. Notice that
\begin{equation}\label{trans}
(\rho_y(p),yx)=\Pi_2^{-1}\circ L_y\circ \Pi_2 (p,x).
\end{equation}
Hence the mapping
$$
(y,p)\lmt \rho_y(p)
$$
is of class $\cc^{k+1}$ from $G\times M$ onto $M$.

To finish we have $\rho_y\circ\rho_{y'}=\rho_{yy'}$ for all $y,y'\in G$, according to Equality (\ref{trans}).

Therefore $(y,p)\lmt \rho_y(p)$ is a transitive and $\cc^{k+1}$ action of the Lie group $G$ on the manifold $M$. Let $H$ be the isotropy group of $p_0$, that is the set of points $x$ of $G$ such that $(p_0,x)$ belongs to $S$, and $G/H$ the manifold of left cosets of $H$. Then 
$$
xH\in G/H\lmt \rho_x(p_0)
$$
is a diffeomorphism from $G/H$ onto $M$, denoted by $\Phi$ in the end of the proof.

\item
It remains to prove that $\Phi$ induces an isomorphism between the Lie algebra of invariant vector fields on $G/H$ and  $\Li$.
Let us denote by $\Pi$ the projection of $G$ onto $G/H$. Then the equality
$$
\Phi\circ \Pi=\Pi_1\circ\Pi_2^{-1}
$$
holds. Recall that $(p_0,e)\in S$. Then $\forall x\in G$
$$
\begin{array}{ll}
\Phi\circ\Pi(x)&=\Phi(xH)=\rho_x(p_0)\\
               &=\Pi_1(\rho_x(p_0),x)\\
               &=\Pi_1\circ\Pi_2^{-1}(x)
\end{array}
$$
Let $Y\in \g$. Then $\Phi_*(\Pi_*Y)=(\Pi_1\circ\Pi_2^{-1})_*Y=L^{-1}(Y)$. This equality has two consequences. The first one is that $\Gamma$ is under $\Phi^{-1}$ equivalent to a set of invariant vector fields on $G/H$. The second one is that all the vector fields of $\Li(\Gamma)$ are complete, since $\Phi$ is a diffeomorphism, and they are related by $\Phi_*$ to complete vector field of $G/H$.
\end{enumerate}

\hfill $\Box$

Thanks to the Orbit Theorem, recalled in Section \ref{defandnot}, we can relax the transitivity assumption, and obtain the following corollary.

\newtheorem{Palais}{Corollary}
\begin{Palais} \label{Palais}
Let $\Gamma$ be a family of $\cc^k$ vector fields on a connected  manifold $M$. If
\begin{enumerate}
\item[(i)] all the vector fields belonging to $\Gamma$ are complete,
\item[(ii)] $\Gamma$ generates a finite dimensional Lie algebra $\Li(\Gamma)$,
\end{enumerate}
then all the vector fields belonging to $\Li(\Gamma)$ are complete,
and the family $\Gamma$ is Lie-determined.
\end{Palais}

\demo
Let $p$ be a point of $M$ and let us denote by $S$ the orbit of
$\Gamma$ through $p$. By the orbit theorem $S$ is a submanifold of
$M$. Moreover every vector field belonging to
$\Gamma$, hence every vector field belonging to $\Li(\Gamma)$, is
tangent to $S$.

Let $\Gamma_S$ stand for the family of restrictions to $S$ of the vector fields of $\Gamma$. Clearly $\Gamma_S$ generates a finite dimensional Lie algebra $\Li(\Gamma_S)$, which is nothing else than the
set of restrictions to $S$ of the vector fields of $\Li(\Gamma)$. By
definition the family $\Gamma_S$ is transitive on $S$ and satisfies the
assumptions of Theorem \ref{tecnic}. Therefore $S$ is diffeomorphic to
a homogeneous space $G/H$, where $G$ is a simply connected Lie Group
whose Lie algebra is isomorphic to $\Li(\Gamma_S)$, and $H$ is a closed
subgroup of $G$. By this diffeomorphism $\Gamma_S$ is related to a
set of invariant vector fields, and $\Li(\Gamma_S)$ to the Lie algebra of invariant vector fields on
$G/H$. Therefore
$$
\forall q\in S \qquad\qquad \rank\Li(\Gamma_S)(q)=\dim G/H=\dim S
$$
and the family $\Gamma_S$ is Lie-determined.

\hfill $\Box$


\section{Application to control systems.}\label{main}

Consider the control affine system 
$$
(\Sigma)\qquad\qquad      \dot{x}=f(x)+\sum_{j=1}^mu_jg_j(x)
$$
where $x$ belongs to the $n$-dimensional connected manifold $M$ and where
$f,g_1,\dots,g_m$ are $\cc^k$ vector fields on $M$, with $k>0$. The control
$u=(u_1,\dots,u_m)$ belongs to $\R^m$.

The family $\Gamma=\{f,g_1,\dots,g_m\}$ is assumed to generate a Lie
algebra denoted by $\Li$.

We also denote by $\Li_0$ the ideal of $\Li$ generated by
$g_1,\dots,g_m$. It is well known that $\Li_0$ is the smallest Lie
subalgebra of $\Li$ containing $g_1,\dots,g_m$ and closed
for the Lie bracket with $f$:
$$
X\in \Li_0 \quad \Lra\quad [f,X]\in \Li_0.
$$

The \textit{dimension} of $\Li_0$ is its dimension as a real Lie
algebra, and its \textit{rank at a point} $p\in M$ is the dimension of the
subspace $\{X(p);\ \ X\in \Li_0\}$ of the tangent space $T_pM$ of $M$
at $p$.
In the particular case where the rank of $\Li_0$ is
constant, it will be refered to as $\rank(\Li_0)$.

If $\Li$ (resp. $\Li_0$) is finite dimensional, then $G$ (resp. $G_0$)
will stand for a (connected and) simply connected Lie group whose Lie algebra is
isomorphic to $\Li$ (resp. $\Li_0$).

\newtheorem{maintheorem}[T1]{Theorem}
\begin{maintheorem} \label{maintheorem}
We assume the family
$\{f,g_1,\dots,g_m\}$ to be transitive.
Then System $(\Sigma)$ is diffeomorphic to a linear system on a Lie group
or a homogeneous space if and only if the vector fields $f,g_1,\dots,g_m$ are complete and generate a 
finite dimensional Lie algebra.

More accurately, under this condition the rank of $\Li_0$ is constant, equal to $\dim (M)$ or $\dim(M)-1$, and:
\begin{enumerate}
\item[(i)] if $\rank(\Li_0)=\dim (M)$, in particular if there exists one
  point $p_0\in M$ such that $f(p_0)=0$, then $(\Sigma)$ is
  diffeomorphic to a linear system on a homogeneous space $G_0/H$ of
  $G_0$;
\item[(ii)] if $\rank(\Li_0)=\dim (M)-1$, then $\Sigma$ is
  diffeomorphic to an invariant system on a homogeneous space $G/H$ of
  $G$.
\end{enumerate}
\end{maintheorem}

\demo Let us prove the sufficiency.

The Lie algebras $\Li$ and $\Li_0$ are finite dimensional and
generated by complete vector fields. By Corollary \ref{Palais} they
are Lie-determined. As $\Li$ is transitive, its rank is
everywhere full. Moreover the rank of $\Li_0$ is constant over $M$,
equal to $\dim(M)$ or $\dim(M)-1$. Indeed $\Li_0$ being Lie-determined, its rank is everywhere equal to the dimension of the zero-time orbit, which is constant, equal to $\dim(M)$ or $\dim(M)-1$ (see \cite{Jurdjevic97}).

\textit{Let us assume $\rank(\Li_0)=\dim (M)$}.
Then we can apply Theorem \ref{tecnic} to the family $\Li_0$. The manifold $M$ is diffeomorphic to a homogeneous space $G_0/H$ of $G_0$, and if we denote by $\Phi$ this diffeomorphism, the tangent mapping $\Phi_*$ induces a Lie algebra isomorphism between  $\Li_0$ and the Lie algebra of invariant vector fields on $G_0/H$. The vector field $\Phi_*f$ satisfies
$$
\forall g\in \Li_0 \qquad [\Phi_*f,\Phi_*g]=\Phi_*[f,g]\in \Phi_*\Li_0.
$$
Since $\Phi_*(\Li_0)$ is equal to the Lie algebra of invariant vector fields on $G_0/H$, and according to Theorem \ref{homogene}, the vector field $\Phi_*f$ is affine, that is $\Phi_*f$ is the projection onto $G_0/H$ of an affine vector field $F$ of $G_0$. This vector field can be chosen to be linear if and only if there is one point $p_0$ in $M$ such that $f(p_0)=0$: in the proof of Theorem \ref{tecnic}, we can choose $p_0$ to be the projection of the identity $e$ of $G_0$. Clearly System $\Sigma$ is diffeomorphic to the linear system
$$
\dot{x}=F(x)+\sum_{j=1}^m u_jY_j(x)
$$
on $G_0/H$, where $F$ stands for $\Phi_*f$, and $Y_j=\Phi_*g_j$ is an invariant vector field for $j=1,\dots, m$.

\textit{We assume now $\rank(\Li_0)=\dim (M)-1$}.
We apply Theorem \ref{tecnic} to $\Li$: the manifold $M$ is diffeomorphic to a homogeneous space $G/H$ of $G$, and under this diffeomorphism $\Li$ is isomorphic to the Lie algebra of invariant vector fields on $G/H$. System $\Sigma$ is obviously diffeomorphic to an invariant system on $G/H$.

\hfill $\Box$

From Theorem \ref{maintheorem} we can deduce the following corollary, stated in the $\cci$ case in order to ensure the existence of $\Li$ and $\Li_0$:

\newtheorem{mainbis}[Palais]{Corollary}
\begin{mainbis}\label{mainbis}
The manifold $M$ is assumed to be $\cci$ and simply connected.

The family $\{f,g_1,\dots,g_m\}$ is assumed to be $\cci$, complete and transitive, and the vector field $f$ to vanish at a point $p_0\in M$.

Then System $\Sigma$ is equivalent to a linear system on a Lie group if and only if
$$
\dim (\Li_0)=\dim (M)
$$
\end{mainbis}

\demo
The necessity part is obvious. Let us prove the sufficient one.

Since $\dim (\Li_0)<\infty$, Theorem \ref{maintheorem} applies, and since $f$ vanishes at one point, we have $\rank (\Li_0)=\dim (M)=\dim (\Li_0)$. Therefore $\Sigma$ is diffeomorphic to a linear system on a homogeneous space $G_0/H$ of $G_0$, with the previous notations. Now the two conditions $\dim (\Li_0)=\dim (M)$ and $M$ simply connected imply $G_0/H=G_0$.

\hfill $\Box$

The assumption that $M$ is simply connected cannot be relaxed. If not, $M$ remains diffeomorphic to a homogeneous space $G/H$, where $H$ is discrete, but $G/H$ is a Lie group if and only if $H$ is normal.


\section{Examples}\label{Exemples}
\subsection{Examples of linear and affine vector fields}
\subsubsection{Inner derivations on matrix Lie groups}
See \cite{Markus81}. Let $G$ be a connected matrix Lie group, that is a connected Lie subgroup of $Gl(n;\R)$, for some $n$. For any matrix $X$ belonging to the tangent space $T_IG$ at the identity $I$, identified to $\g$, the mapping $M\lmt XM$ defines a right invariant vector field. We can also associate to the inner derivation $D=-\ad X$ the linear vector field $\X$ defined by
$$
\X(M)=XM-MX
$$
Hence in the inner derivation case, a linear system on $G$ writes
$$
\dot{M}=XM-MX+\sum_{j=1}^mu_jY_jM
$$
where $M\in G$, and $X,Y_1,\dots,Y_m\in \g$.


\subsubsection{Affine vector fields on the sphere $S^n$, $n\geq 2$}

The sphere $S^n$ is diffeomorphic to the homogeneous space $SO_{n+1}/SO_n$, where $SO_n$ is identified with the closed subgroup
$$
H=\left\{
\left(
\begin{array}{c|c}
 1&0\ \ \hdots \ \ 0\\
\hline\\
0&\\
\vdots& N\\
0&
\end{array}
\right);\ N\in SO_n
\right\}
$$
of $SO_{n+1}$.

On the one hand the Lie algebra $\mathfrak{so}_{n+1}$ is semi simple, since $n+1\geq 3$, so all its derivations are inner. 

On the other hand the subgroup $H$ is connected, and, following Proposition \ref{Projecthomogene}, a linear vector field $\X$ on $SO_{n+1}$ is related to a vector field on $SO_{n+1}/H$ if and only if its Lie algebra $\h$ is invariant under $D=-\ad (\X)$.

But $\X$ being of the form
$$
\X(M)=XM-MX, \qquad\qquad M\in SO_{n+1}
$$
for some $X\in \mathfrak{so}_{n+1}$, this condition turns out to
$$
\forall Y\in \h \qquad [X,Y]\in \h
$$
and a straightforward computation shows that it holds if and only if $X\in\h$.

The flow of $\X$ is given by
$$
\varphi_t(M)=e^{tX}Me^{-tX}
$$
and its projection onto $SO_{n+1}/H$ is equal to
$$
e^{tX}Me^{-tX}H=e^{tX}MH
$$
since $\forall t\in\R$, $e^{-tX}\in H$. This proves that the vector fields $\X$ and $X$ have the same projection on $S^n\sim SO_{n+1}/H$.

In conclusion the only affine vector fields on the sphere $S^n$ are the invariant ones. They are the vector fields defined by
$$
f(x)=Ax, \qquad x\in S^n
$$
where $A\in \mathfrak{so}_{n+1}$, that is $A'=-A$.

\noindent{\bf Remark}. The group $SO_{n+1}$ is not simply connected but as it is semi simple, this is not a restriction. Indeed any derivation $D$ is inner, hence the associated linear vector field always exists, despite the lack of simply connectedness.


\subsubsection{The general inner derivation case}
The previous phenomena is due to the fact that the normalizer of $\h\sim\mathfrak{so}_n$ in $\mathfrak{so}_{n+1}$ is itself. To see this let us consider the general inner derivation case on a connected Lie group $G$.

Let $D=-\ad X$, with $X\in \g$, be an inner derivation. Recall from Section \ref{Lin} that the linear vector field associated to $D$ is equal to
$$
\X=X+\mathcal{I}_*X.
$$
Let $H$ be a closed and connected subgroup of $G$. Then, following again Proposition \ref{Projecthomogene}, the linear vector field $\X$ is related to a vector field on $G/H$ if and only if its Lie algebra $\h$ is invariant under $D=-\ad (\X)$. But
$$
\begin{array}{ll}
& \forall Y\in\h \qquad [Y,\X]=[Y,X]\in H\\
\Longleftrightarrow & X\in \mbox{norm}_{\g}\h
\end{array}
$$
Therefore $\X$ can be projected on $G/H$ if and only if $X$ belongs to the normalizer of $\h$ in $\g$.

\vskip 0.2cm

An example of linear vector field whose derivation is not inner is given in the next section.

\subsection{Example of equivalence}

Consider the system in $\R^2$
$$
\Sigma=\left\{
\begin{array}{ll}
\dot{x} & =y^2\\
\dot{y} & =u
\end{array}
\right.
$$
Let
$$
f=y^2\ddx{x}, \qquad g_1=\ddx{y}, \qquad g_2=[g_1,f]=2y\ddx{x},\  \mbox{ and
} g_3=[g_1,g_2]=2\ddx{x}.
$$
Clearly the system can be written
$$
\dot{p}=f(p)+ug(p) \qquad\qquad\mbox{where }\qquad p=(x,y).
$$
The brackets $[g_2,g_3]$, $[f,g_2]$, and $[f,g_3]$ vanishes, and with
the notations of Section \ref{main} we have
$$
\Li=Sp\{f,g_1,g_2,g_3\} \ \mbox{ and } \ \Lo=Sp\{g_1,g_2,g_3\}.
$$

All these vector fields are complete, $\Lo$ is isomorphic to the
Heisenberg Lie algebra, and its rank is everywhere full, therefore
Theorem \ref{maintheorem} applies: System $\Sigma$ is equivalent to a linear
system on a homogeneous space of the Heisenberg group.

Let us compute this equivalence. The Heisenberg group is
$$
G=\left\{\begin{pmatrix}1&y&z\\0&1&x\\0&0&1\end{pmatrix};\ \ (x,y,z)\in\R\right\}
$$
and its Lie algebra $\g$ is spanned by the right invariant
vector fields
$$
X=\begin{pmatrix}0&0&0\\0&0&1\\0&0&0\end{pmatrix},\quad
Y=\begin{pmatrix}0&1&x\\0&0&0\\0&0&0\end{pmatrix},
\quad \mbox{and }
Z=\begin{pmatrix}0&0&1\\0&0&0\\0&0&0\end{pmatrix},
$$
that can be written in the canonical coordinates
$$
X=\ddx{x}, \qquad
Y=\ddx{y}+x\ddx{z}, \qquad \mbox{and }\ 
Z=\ddx{z}.
$$

The derivation $D$ on $\g$ should verify $DX=Y$ and
$DY=DZ=0$.
Let $\X$ be the vector field on $G$ defined by
$$
\X=x\ddx{y}+\frac{1}{2}x^2\ddx{z}
$$
It is easy to see that $\X$ is linear and that $-\ad(\X)=D$. The
system $\dot{q}=\X(q)+uX(q)$ can be written in $\R^3$
$$
\left\{
\begin{array}{ll}
\dot{x} & =u\\
\dot{y} & =x\\
\dot{z} & =\frac{1}{2}x^2
\end{array}
\right.
$$
We are looking for a subgroup $H$ of $G$ for
which $\Sigma$ is equivalent to a linear system on $G/H$. The algebra
of vector fields of $\Lo$ that vanishes at $(0,0)$ is
spanned by $g_2$, thus the Lie algebra of $H$ should be spanned by $Y$.

Let $H$ be the closed, but not normal, subgroup of $G$
$$
H=\left\{\begin{pmatrix}1&y&0\\0&1&0\\0&0&1\end{pmatrix};\ \ y\in\R\right\}.
$$
The projection of $G$ onto $G/H$ is equivalent to the projection
$(x,y,z)\lmt(x,z)$ from $\R^3$ onto $\R^2$, and the linear system on
$G/H$ is therefore equivalent to
$$
\Sigma'=\left\{
\begin{array}{ll}
\dot{x} & =u\\
\dot{z} & =\frac{1}{2}x^2
\end{array}
\right.
$$
To finish, $\Sigma'$ is equivalent to $\Sigma$ under the linear transformation of $\R^2$, $(x,z)\lmt (2z,x)$.

\subsection{Generalization of the previous example}

Let $P(y)$ be a polynomial. Then the system in $\R^2$
$$
\left\{
\begin{array}{ll}
\dot{x} & =P(y)\\
\dot{y} & =u
\end{array}
\right.
$$
satisfies the assumptions of Theorem \ref{maintheorem}: let
$$
f=P(y)\ddx{x}, \qquad g_1=\ddx{y}, \qquad g_2=[g_1,f]=P'(y)\ddx{x}
$$
and, by induction:
$$
g_{k+1}=[g_1,g_k]=P^{(k)}(y)\ddx{x}.
$$
The vector field $g_k$ vanishes as soon as $k>\deg(P)+1$, and so does the other brackets. Therefore 
$$
\Li=Sp\{f,g_1,\dots,g_{\deg(P)+1}\} \ \mbox{ and } \ \Lo=Sp\{g_1,\dots,g_{\deg(P)+1}\}.
$$
All these vector fields are complete, the system satisfies clearly the rank condition, therefore Theorem \ref{maintheorem} applies.

\hfill $\Box$

\vskip 0.5cm

\noindent{\bf Acknowledgments.} The author wishes to express his thanks to Professor Witold Respondek for drawing his attention to equivalence problems, and for stimulating conversations.


\end{document}